\documentclass{amsart}
\usepackage{amsfonts,amssymb,amscd,amsmath,enumerate,verbatim,calc}
\usepackage{tikz-cd}
\usepackage[all]{xy}
 \usepackage{hyperref}

\newcommand{\xhookrightarrow}{\lhook\joinrel\longrightarrow}

\newcommand{\rt}{\rightarrow}

\newcommand{\Der}{\operatorname{Der}}
\newcommand{\rank}{\operatorname{rank}}

\newcommand{\Hom}{\operatorname{Hom}}

\newcommand{\Tor}{\operatorname{Tor}}
\newcommand{\Supp}{\operatorname{Supp}}

\theoremstyle{plain}

\newtheorem{theorem}{Theorem}[section]
\newtheorem{corollary}[theorem]{Corollary}

\theoremstyle{definition}
\newtheorem{definition}[theorem]{Definition}
\newtheorem{s}[theorem]{}
\newtheorem{remark}[theorem]{Remark}
\newtheorem{example}[theorem]{Example}

\theoremstyle{remark}

\begin{document}

\title{ On some Rings of differentiable type  }
 \author{Sayed Sadiqul Islam }
  \email{22d0786@iitb.ac.in}
 \author{Tony J. Puthenpurakal}
\email{tputhen@math.iitb.ac.in}
\date{\today}
\address{Department of mathematics, Indian Institute of Technology Bombay, Powai, Mumbai 400076, India}

\date{\today}

\begin{abstract}
Let $K$ be a field of characteristic 0 and $S=K[x_1,\ldots,x_m]/I$ be an affine domain. Consider $R=S_P$ where $P\in Spec(S)$ such that $R$ is regular. In this paper we construct a field $F$ which is contained in  $R$ such that
\begin{enumerate}
  \item The residue field of $R$ is a finite extension of $F$.
  \item $D_F(R)$,  the ring of $F$-linear differential operators on $R$ is left and right Noetherian with finite global dimension.
  \item  The Bernstein class of $D_F(R)$ is closed under localization at one element of $R$.
\end{enumerate}
 We also prove a similar  result for $R^h$, the Henselization of $R$. As an application we prove that $\frac{D_F(R)}{D_F(R)P}\cong E(\kappa(P))$ where $E(\kappa(P))$ is the injective hull of the residue field of $R$.
\end{abstract}

\maketitle

\section{Introduction}
Let $K$ be a field of characteristic zero, and $R$ be a $K$-algebra. let $D(R)=D_K(R)$ be the ring of $K$-linear differential operators on $R$. There has been a lot of research on the finiteness properties of rings of differential operators. The following are example of $R$ related to $D(R)$ as a left and right Noetherian with finite global dimension:
\begin{enumerate}[\rm (1)]
\item Polynomial rings over $K$ in $n$ variables. In this case $D(R)=A_n(K)$, the $n^{th}-$ Weyl algebra over $K$. We have that global dimension of $D(R)$ is equal to $n$; see \cite[Chapter 2, Theorem 3.15]{B}.

\item Power series rings over $K$ in $n$ variables. In this case global dimension is $n$; see \cite[Chapter 3, Proposition 1.8]{B}.
\item The ring of convergent power series over $\mathbb{C}$ in $n$ variables; see \cite[Page 197]{B}.
\item Let $(V,\pi V,K)$ be a DVR of mixed characteristic $p>0$, and Let $F$ denote its fraction field. Let $S=V[[x_1,\ldots,x_{n+1}]]\otimes_V F$ be the tensor product of the power series ring
with coefficients in $V$ and $F$. Let $R=S/(f)S$ be a regular ring where $f=\pi-h$ for an element $h$ in the square of maximal ideal of $V[[x_1,\ldots,x_{n+1}]]$. In this case also $D(R)$ is Noetherian with finite global dimension; see \cite[Page 4]{NB}.
\item Localization of complete local domains of equal characteristic zero at certain elements; see \cite{Tony}.
\end{enumerate}
 However $D(R)$ is not always Noetherian. Consider in the space $\mathbb C^3$ with the coordinates $x_1,x_2,x_3$ the surface $X$ defined by the equation $x^3_1+x^3_2+x^3_3=0$, then the ring of differential operator on $X$ is not Noetherian. This example is due to Bernstein- Gelfand\cite[Theorem 1]{Ber}. In this article we provide a new infinite class of Noetherian domains $R$ with $D(R)$ is both left and right Noetherian and with finite global dimension.

 We now state the following results of this paper.
\s \label{Hypothesis-2} \textit{Setup:}
\label{Exm}
Let $K$ be a field of characteristic 0 and $S=K[x_1,\ldots,x_m]/I$ be an affine domain. Consider $R=S_P$ where $P\in Spec(S)$ such that $R$ is regular. We shall show that $R$ contains a field $F$ with $\kappa(P)=S_P/PS_P$ a finite extension of $F$. Let $\Der_F(R)$ be the set of $F$-linear derivations of $R$ and $D_F(R)$ be the ring of $F$-linear differential operators on $R$. Our result is as follows:
\begin{theorem}
\label{Main-1.1}
    Let $R$ and $F$ be as in Setup \ref{Exm}. Then,
    \begin{enumerate}[\rm (1)]
        \item The residue field of $R$ is a finite extension of $F$.
        \item $\Der_F(R)$ is a finitely generated projective $R$-module of $\rank$ equals to $\dim R$.
        \item The ring $D_F(R)$ is a left and right Noetherian and the Bernstein class of $D_F(R)$ is closed under localization at one element of $R$.
    \end{enumerate}
\end{theorem}
\s \label{Hypothesis-3} \textit{Setup:}
Let $R$ and $S$ be as in Setup \ref{Exm}. Denote by $R^h$, the Henselization of $R$. Since $R\hookrightarrow R^h$, the field $F$ in \ref{Exm}
 is contained in $R^h$. We prove a similar result.

 \begin{theorem}
 \label{Main-2}
 For the setup \ref{Hypothesis-3}, we have the following;
\begin{enumerate}[\rm (1)]
    \item The residue field of $R^h$ is a finite extension of $F$.
\item $\Der_F(R^h)$ is a finitely generated projective $R^h$-module of rank same as dimension of $R^h$.
\item The ring $D_F(R^h)$ is a left and right Noetherian and the Bernstein class of $D_F(R^h)$ is closed under localization at one element of $R^h$.
\end{enumerate}
 \end{theorem}
 We also prove the following result as an application of \ref{Main-1.1};
 \begin{theorem}
 \label{1.5}
      Let $R$ and $F$ be as in Setup \ref{Exm}. Then,
      \begin{enumerate} [\rm (1)]
      \item $H_{P}^{\dim R}(R)\cong E(\kappa(P))$ is a $D_F(R)$-module in the Bernstein class.
          \item $\frac{D_F(R)}{D_F(R)P}\cong E(\kappa(P))$,
       where $E(\kappa(P))$ is the injective hull of the residue field of $R$.
      \end{enumerate}

 \end{theorem}
Our paper comprises six sections. Section 2 provides an overview of $D$-modules. In section 3, we present crucial theorems essential for proving our results. Sections 4, 5, 6 and 7 are where we explain our proofs.
 \section{Background on D-modules}

Let $R$ be a commutative $A$-algebra. The ring of $A$-linear differential operators of $R$ is a subring $D_A(R)\subseteq \Hom_{A}(R,R)$ whose elements are defined inductively as follows: the differential operators of order zero are defined by the multiplication by elements of $R$, i.e. $D^0(R)\cong R$. Suppose we have defined operators of order $< n$. An operator $\delta \in \Hom_{A}(R,R)$ is of order less than or equal to $n$ if $[\delta,r]=\delta r-r\delta$ is an operator of order less than or equal to $n-1$ for all $r\in R$. We have a filtration $D_A^0(R)\subset D_A^1(R)\subset\ldots$ and ring of differential operators is defined as
$$D_A(R)=\bigcup_{m\geq 0} D_A^m(R).$$
\begin{example}
    If $R$ is either a polynomial ring $A[x_1,\ldots,x_n]$ or the formal power series ring $A[[x_1,\ldots,x_n]] $ with coefficients in $A$, then ring of $A$-linear differential operators is :
    $$D_A(R)=R\left\langle\left.\frac{1}{t!} \frac{d^t}{d x_i^t} \right\rvert\, i=1, \ldots, n ; t \in \mathbb{N}\right\rangle.$$
\end{example}
\s \label{2.2}  Let $M$ be a $D_A(R)$-module and $f\in R$. Then $M_f$ is also a $D_A(R)$-module, where the action of a differential operator $\partial\in D_A(R)$ on $\frac{v}{f^t}\in M_f$ is defined inductively as follows: If $\partial\in D^0_A(R)$ has order 0 then $\partial\cdot\frac{v}{f^t}=\frac{\partial\cdot v}{f^t}$. Now suppose that action of every element in $D^n_A(R)$ has been defined. Let $\partial\in D^{n+1}_A(R) $. Then,
 $$\partial\cdot \frac{v}{f^t}=\frac{\partial\cdot v-[\partial,f^t]\cdot \frac{v}{f^t}}{f^t} $$
 With this $D_A(R)$-module structure on $M_f$, the map $M\rightarrow M_f$ is a morphism of $D_A(R)$-modules.

\s \label{2.3}Let $I\subseteq R$ be an ideal of $R$ generated by $\underline f= f_1,f_2,\ldots,f_r\in R$, and let $M$ be an $R$-module. Then the \v Cech complex of $M $ with respect to  $\underline f$ is defined by
 $$
  \Check{C}^{\bullet}(\underline f,M):\  0\rightarrow M\rightarrow \bigoplus_iM_{f_i}\rightarrow \bigoplus_{i,j}M_{f_if_j}\rightarrow\ldots \rightarrow M_{f_1\ldots f_r}\rightarrow 0
 $$
 where the maps on every summand are localization map upto a sign. The local cohomology module of $M$ with support on $I$  is defined by
 $$H_I^i(M)=H^i( \Check{C}^{\bullet}(\underline f,M))$$
 It follows from this construction that every local cohomology module over a $D_A(R)$-module is again a $D_A(R)$-module.
 \begin{definition}
        We say that a associative ring $R$ is filtered if there exists an ascending filtration $\Sigma_0\subset\Sigma_1\subset\Sigma_2\ldots$ of additive subgroups such that $1\in \Sigma_0$, $\bigcup \Sigma_i=R$ and $\Sigma_i\Sigma_j\subset\Sigma_{i+j}$ for every $i$ and $j$. We denote by $gr^{\Sigma}(R)$ the associated graded ring $\Sigma_0\oplus \Sigma_1/\Sigma_0\oplus \Sigma_2/\Sigma_1\oplus\ldots$
    \end{definition}
    We have the following definition of a ring of differentiable type from \cite{Mar};
    \begin{definition}
     A filtered ring $R$ is a ring of differentiable type if its associated graded ring is commutative Noetherian regular with unity and pure graded dimension.
    \end{definition}
\s    Let $R$ be a filtered ring of differentiable type with filtration $\{\Sigma_i\}$. A filtration $\{\Gamma_i\}$ compatible with $\{\Sigma_i\}$ on a left $R$-module $M$ is called good if the associated graded $gr^{\Sigma}(R)$-module $gr^{\Gamma}(M)$ is finitely generated. It is well known that a left $R$-module $M$ can be equipped with good filtration iff it is finitely generated. If $M$ is a finitely generated left $R$-module we can choose some good filtration $\Gamma$ on $M$. Now the associated graded $gr^{\Sigma}(R)$-module $gr^{\Gamma}(M)$ is finitely generated. Then we define:
    $$d_{\Gamma}(M)=\dim \Supp(gr^{\Gamma}(M))$$
    In general, $M$ can be equipped with many good filtrations. But it turns out that the integer $d_{\Gamma}(M)$ do not depend on the choice of $\Gamma$, i.e. we have the equalities $d_{\Gamma}(M)=d_{\Gamma^{\prime}}(M)$ for all pairs of good filtrations $\Gamma$ and $\Gamma^{\prime}$ on $M$. This common integer is therefore denoted by $d(M)$ and called the dimension of $R$-module $M$.
    \begin{definition}
    A ring $R$ has finite weak global dimension $w$ if $\Tor^R_v(N,M)=0$ for all $v>w$ and every right $R$-module $N$ and left $R$-module $M$, while $\Tor^R_w(N,M)\neq 0$ for some pair.
\end{definition}

Let $R$ be a ring of differentiable type with pure graded dimension $\mu$ and assume that weak global dimension of $R$ is $w$. We recall the following corollary from \cite[Chapter 2, 7.2]{B},
\begin{corollary}
    $d(M)\geq \mu-w$ for every non-zero finitely generated R-module $M$.
\end{corollary}
\begin{definition}
    The Bernstein class consists of all finitely generated $R$-modules $M$ for which equality holds in the above corollary i.e. $d(M)= \mu-w$ (or if $M = 0$).
\end{definition}
 \section{Preliminaries}
If $A$ is a $k$-algebra, $\Omega_{A/k}$ denotes the module of  K\"{a}hler differentials and let the universal derivation be $d_{A/k}\colon A \rt \Omega_{A/k} $. For $A$-module M we have that $\Der_k(A,M)\cong \Hom_A(\Omega_{A/k},M)$.

The following two results are useful:
\begin{theorem}(First fundamental exact sequence)\cite[Theorem 25.1]{Mat}
\label{Fun}
    A composite $k\xrightarrow{f} A\xrightarrow{g} B$ of ring homomorphisms leads to an exact sequence of $B$-modules
   $$ \Omega_{A/k}\otimes_{A} B\xrightarrow{\alpha} \Omega_{B/k}\xrightarrow{\beta} \Omega_{B/A}\xrightarrow{} 0$$

   where the maps are given by $\alpha(d_{A/k}a\otimes b)=b\ d_{B/k}g(a)$ and $\beta(d_{B/k}b)=d_{B/A} b$ for $a\in A$ and $b\in B$.
\end{theorem}
\begin{theorem} (Second fundamental exact sequence)\cite[Theorem 25.2]{Mat}
    \label{2nd fun}
  Consider $k\xrightarrow{f} A\xrightarrow{g} B$ when $g$ is surjective; set $\ker g=\mathfrak m$, so $B=A/\mathfrak m$. Then we have an exact sequence as follows:

  $$\mathfrak m/\mathfrak m^2\xrightarrow{\delta} \Omega_{A/k}\otimes_A B\xrightarrow{\alpha}\Omega_{B/k}\rightarrow 0$$
  where $\delta$ is the $B$-linear map defined by $\delta(\Bar{x})=d_{A/k}x\otimes 1$.
\end{theorem}
In general $\Omega_{A/F}$ need not be finitely generated $A$-module. However we have the following result for $\Der_F(A)$.
\begin{theorem}\cite[Theorem 30.7]{Mat}
\label{Tra}
    $(R,\mathfrak m)$ be an $n$-dimensional Noetherian local integral domain containing $\mathbb Q$, and $F\subset R$ be a subfield such that $tr.deg_F(R/\mathfrak m)=r<\infty$. Then $\Der_F(R)$ is isomorphic to a submodule of $R^{n+r}$, and is therefore a finite $R$-module, with
    $$\rank Der_F(R)\leq \operatorname{dim} R+tr.deg_F(R/\mathfrak m)$$
    \end{theorem}
    The following result gives a sufficient condition for $\Der_F(R)$ to be free.
\begin{theorem}\cite[Theorem 30.6]{Mat}
\label{Der}
    Let $(R,\mathfrak m,F)$ be a regular local commutative Noetherian ring with unity of dimension $n$ containing a field $F_0$. Suppose that $F$ is an algebraic separable extension $F_0$. Let $\widehat{R}$ denote the $\mathfrak m$-adic completion of $R$. Let $x_1,\ldots,x_n$ be a regular system of parameters of $R$. Then, $\widehat{R}=F[[x_1,\ldots,x_n]]$ is the power series ring with coefficients in $F$, and $\Der_F\widehat{R}$ is a free $\widehat{R}$-module with basis $\partial/\partial x_i$ for $1\leq i\leq n$. Moreover, the following conditions are equivalent:
    \begin{enumerate}[\rm (1)]
        \item $\partial/\partial x_i$ maps $R$ into $R$ for $1\leq i\leq n$, so that they can be considered as elements of $\Der_{F_0}(R)$;
        \item $\Der_{F_0}(R)$ is a free $R$-module of $\rank$ $n$;
        \item $\rank$ $\Der_{F_0}(R)=n$.
    \end{enumerate}
    \end{theorem}
    \begin{corollary}\cite[Corollary 1.2.2.3]{Mar 2}
    \label{3.5}
        Under the hypotheses above, assume that any of the equivalent conditions in theorem \ref{Der} are satisfied. Then the following holds:
        $$\widehat{R}\otimes_RD_{F_0}(R)\cong D_{F_0}(\widehat{R})\cong D_{F_0}(R)\otimes_R\widehat{R}$$
    \end{corollary}
    The following result is due to Lyubeznik;
    \begin{theorem}\cite[Proposition 2.3]{Ly}
    \label{3.6}
        Let $K$ be a field of characteristic 0, let $R=K[[x_1,\ldots,x_n]]$ be the ring of formal power series ring in $n$ variables over $K$ and let $\mathfrak m$ be maximal ideal of $R$. Then, as an $R$-module $\frac{D_K(R)}{D_K(R)\mathfrak m}$ is isomorphic to $E_R(K)$, the injective hull of the residue field of $R$.
    \end{theorem}
    \begin{remark}
       \label{rem}

       Let $R$ be an integral domain and $L$ be quotient field of $R$.

       \begin{enumerate}

            \item   Let $M$ be $R$-module of $\rank$ $r$. Then $M\otimes L\cong L^r$. But $\Hom_R(M,R)\otimes L\cong \Hom_L(L^r,L)\cong L^r$. Hence $\Hom_R(M,R)$ and $M $ have the same $\rank$.
       \item  Let $M$ is $R$-module and $P\in Spec(R)$. We note that $R$ and $R_P$ have same quotient field, denoted as $L$. Since   $M_P\otimes_{R_P} L\cong M\otimes_R R_P\otimes_{R_P} L\cong M\otimes_R L$, so $M$ and $M_P$ have the same $\rank$.
     \end{enumerate}
    \end{remark}
\begin{definition}
    A local homomorphism $\phi:(R,\mathfrak m,k)\rightarrow (S,\mathfrak n,l)$ of local rings is unramified provided $S$ is essentially of finite type over $R$(that is, $S$ is localization of some finitely generated $R$-algebra) and the following properties hold.
    \begin{enumerate}[\rm (1)]
        \item $\mathfrak mS=\mathfrak n$, and
        \item $S/\mathfrak mS$ is a finite separable field extension of $R/\mathfrak m$.
    \end{enumerate}
    If, in addition, $\phi$ is flat, then we say $\phi$ is \'{e}tale. (We say also that $S$ is an unramified, respectively, \'{e}tale extension of $R$.) A pointed \'{e}tale neighborhood is an \'{e}tale extension  $(R,\mathfrak m,k)\rightarrow (S,\mathfrak n,l)$ inducing an isomorphism on residue fields.
    \end{definition}

    The isomorphism classes of pointed \'{e}tale neighborhoods form a direct system. The reason for this is that if  $R\rightarrow S$ and $R\rightarrow T$ are pointed \'{e}tale neighborhoods then there is at most one homomorphism $S\rightarrow T$ making the obvious diagram commute. We refer \cite[Chapter 3]{Ive} for details.
    \begin{definition}
        The Henselization $R^h$ of $R$ is the direct limit of a set of representatives of isomorphism classes of pointed \'{e}tale neighborhoods of $R$.
    \end{definition}
    \begin{theorem}\cite[Theorem 43.9]{Nag}
    \label{Naga}
        Let $(R,\mathfrak m)$ be a local ring and let $R^h$ be its Henselization. For a finite number of elements $b_i$ of $R^h$, there is an element $a$ of $\mathfrak mR^h$ which is a root of a monic polynomial $f(x)=x^r+c_1x^{r-1}+\ldots+c_r$ over $R$ with $c_{r-1}\notin \mathfrak m,\ c_r\in \mathfrak m$, such that the ring $R'=R[a]_{(\mathfrak m+aR[a])}$ contains $b_i$. Furthermore  $R^h$ contains $R'$.
    \end{theorem}

    We recall the following two results from \cite{NB}.
    \begin{theorem}
    \label{BT}
        Let $R$ be a regular commutative Noetherian ring with unity that contains a field $F$ of characteristic 0 satisfying the following conditions:
        \begin{enumerate}[\rm (1)]
            \item R is equidimensional of dimension $n$;
            \item  every residual field with respect to a maximal ideal is an algebraic extension of $F$;
            \item $\Der_F(R)$ is a finitely generated projective $R$-module of rank $n$ such that $R_\mathfrak m\otimes_R\Der_F(R)=Der_F(R_\mathfrak m)$ for each maximal ideal $\mathfrak m$ of $R$.
        \end{enumerate}
    Then, the ring of $F$-linear differential operators $D_F(R)$ is a ring of differentiable type of weak global dimension equal to $\dim(R)$. Moreover, the Bernstein class of $D_F(R)$ is closed under localization at one element of $R$.
    \end{theorem}
    Let $D_F(R)$ be as in the last theorem. Then,
\begin{theorem}
\label{3.12}
    If M is a left $D_F(R)$-module in the Bernstein class. Then, each local cohomology module has a natural structure of $D_F(R)$-module and it belongs to Bernstein class of $D_F(R)$-module.
\end{theorem}
\section{Existence of the field of  \ref{Hypothesis-2}}

\label{Fi}
Let $R$ be as in setup \ref{Hypothesis-2}. We first show the existence of a field $F$ contained in $R$ with $\kappa(P)$ a finite extension of $F$. As $S/P$ is a finitely generated $K$-algebra, there exist variables $y_1,\ldots,y_r$ such that $T=K[y_1,\ldots,y_r]\xhookrightarrow{} S/P$ and $S/P$ is integral over $K[y_1,\ldots,y_r]$ (By Noether normalization). So quotient field of $S/P$ is a finite extension of quotient field of $K[y_1,\ldots,y_r]$. Set $F$ to be the quotient field of $K[y_1,\ldots y_r]$. Therefore $S_P/PS_P$ is a finite extension of $F$. We claim that $F$ is contained in $R=S_P$. Let $\phi$ denotes the inclusion map between $T$ and $S/P$. If $\phi(y_i)=\overline{a_i}$ then the map $\psi:T\rightarrow S$ given by $\psi(y_i)=a_i$ is injective. Composition of $\psi$ with the natural  map from $S$ to $S_P$ is injective. Let this map be $h:T\rightarrow S_P$.
\begin{equation*}
    \begin{tikzcd}
        T \arrow{d}{} \arrow{r}{h} & S_p \\
        T_{(0)} \arrow[dashed]{ru}{g}
    \end{tikzcd}
\end{equation*}

Every non zero element in $T$ is non zero in $S/P$ and hence unit in $S_P$. So by universal property of localization there exists a map $g:T_{(0)}\rightarrow S_P$ such that the above diagram is commutative. Hence $T_{(0)}=F$ is contained in $S_P$ as $g$ is injective.
\section{Proof of theorem \ref{Main-1.1}}
\label{DerR}
From the discussions above we can see that $S_P/PS_P$ is finite extension of $F$ and this proves first part of \ref{Main-1.1}.
To prove second part of \ref{Main-1.1}, it is enough to prove that $\Der_F(R)$ is free of $\rank$ $n = \dim R$ because $R$ is local ring.
\label{Free}
\begin{proof}
If $S$ is not regular, let $Sing(S)$, the singular locus of $S$, be defined by an ideal $J$ in $S$. Note that $J\neq 0$ and $J\not \subset P$. So there exists $  \alpha\in J$ such that $\alpha\notin P$. Then $S_\alpha=V$ is a regular affine domain and $R=S_P=V_{PV}$. From this observation we may assume that $S$ is regular affine domain.

It suffices to show $\rank \Der_F(R)=n$ by \ref{Der}. Using \ref{Tra}, $\rank \Der_F(R)\leq n$. We have that $T=K[y_1,\ldots,y_r]\subset S$ and $K\hookrightarrow T\hookrightarrow S$. By \ref{Fun}
 we get the following exact sequence of $S$-modules
 $$\Omega_{T/K}\otimes_T S\rightarrow \Omega_{S/K}\rightarrow \Omega_{S/T}\rightarrow 0 $$
  $\Omega_{T/K}= T^r$ and $\Omega_{S/K}$ is a projective $S$-module of $\rank$ equals to $\dim S=m$. So the exact sequence becomes
 $$S^r=T^r\otimes_T S\rightarrow \Omega_{S/K}\rightarrow \Omega_{S/T}\rightarrow 0$$
 Let $L$ be quotient field of $S$ and $t=\rank \Omega_{S/T}$. Tensoring the above sequence with $L$ we get the exact sequence
 $$L^r\rightarrow L^m\rightarrow L^t\rightarrow 0$$
 Hence $t\geq m-r=\dim S- \dim S/P=\operatorname{ht} P=\dim R=n$.
 Hence $\rank\ \Omega_{S/T}\geq n$.
 Since $\Der_T(S)\cong \Hom_S(\Omega_{S/T},S)$, remark \ref{rem} implies that $\rank \Der_T(S)=\rank\Omega_{S/T}\geq n$. Therefore using \ref{rem} we get $\rank_R \Der_T(S)_P\geq n$. This implies that $\rank  \Der_T(S_P)\geq n$ as $\Der_T(S)_P\subset \Der_T(S_P)$. Clearly $\Der_F(S_P)\subset \Der_T(S_P)$. We claim that this is an equality. Let $D\in \Der_T(S_P)$ and $\frac{a}{s}\in F,\  u\in S_P$. If
 $$v=\frac{a}{s}\ u\implies s v= a u\implies D(sv)=D(au)\implies sD(v)=aD(u)$$
 Hence we get that $D(v)=\frac{a}{s} D(u)$. Thus $D\in \Der_F(S_P)$. So $\Der_F(S_P)= \Der_T(S_P) $. Therefore $\rank \Der_F(R)\geq n$. Hence $\rank \Der_F(R)= n$ and this proves that $\Der_F(R)$ is a free $R$ module by \ref{Der}.
 \end{proof}

     As a corollary of Theorem \ref{BT} we get $3$ of \ref{Main-1.1} that is $D(R)$ is a left and right Noetherian with global dimension $n$. Moreover, the Bernstein class of $D(R)$ is closed under localization at one element of $R$.

\section{Proof of theorem \ref{Main-2}}
 Since the residue fields of $R^h$ and $R$ are same, the residue field of $R^h$ is finite extension of $F$. Only thing is to prove that $\Der_F(R^h)$ is a finitely generated free $R^h$-module of rank same as $\dim R^h=\dim R=n$.
     \begin{proof}
Using \ref{Naga}, if $a\in R^h$, then $a\in \frac{R[x]_{(\mathfrak m_R,x)}}{(f(x))}=T\subseteq R^h$ for some $f(x)=x^n+a_{n-1}x^{n-1}+\ldots+a_1x+a_0$ with $a_0\in\mathfrak m$, $a_1\notin \mathfrak m$.

Let $L=R[x]_{(\mathfrak m_R,x)}=S[x]_{(P,x)}$ and $T=\frac{R[x]_{(\mathfrak m,x)}}{(f(x))}=L/(f)$. Since $F\hookrightarrow R$, we have that $F\hookrightarrow L$ and $L/\mathfrak m_L=R/\mathfrak m$ is a finite extension of $F$.

Using \ref{2nd fun} for $F\rightarrow L\rightarrow L/(f)$, we get the exact sequence
$$(f)/(f)^2\rightarrow \Omega_{L/F}\otimes_L T\rightarrow \Omega_{T/F}\rightarrow 0$$

Applying $\Hom_T(\_,T)$, the exact sequence becomes
$$0\rightarrow \Hom_T(\Omega_{T/F},T)\rightarrow \Hom_T(\Omega_{L/F}\otimes_L T,T)\rightarrow \Hom_T((f)/(f)^2,T)$$
Since  $\Hom_T(\Omega_{T/F},T)\cong \Der_F(T,T)$, $\Hom_T(\Omega_{L/F}\otimes_L T,T)\cong \Hom_L(\Omega_{L/F},T)\cong \Der_F(L,T)$ and $(f)/(f)^2\cong T$, we have the exact sequence
$$0\rightarrow \Der_F(T,T)\rightarrow \Der_F(L,T)\rightarrow \Hom_T(T,T)=T$$
 We claim that $\rank_T \Der_F(L,T)\geq {n+1}.$ Assume the claim for the moment. The last exact sequence implies $\rank_T \Der_F(T)\geq n$. Hence $\rank \Der_F(T)=n$ as \ref{Tra} implies $\rank \Der_F(T)\leq n$. If $a\in T$, then $\frac{\partial}{\partial x_i}(a)\in T$ using \ref{Der} for all $i$. Therefore, $\frac{\partial}{\partial x_i}(a)\in R^h$ because $T\subseteq R^h$. Since this is true for all $a\in R^h$, we get  $\frac{\partial}{\partial x_i}\in \Der_F(R^h)$ for each $i$. Again using \ref{Der}, $\Der_F(R^h)$ is free $R^h$-module of rank $n$.

 Now we prove our claim that $\rank_T \Der_F(L,T)\geq {n+1}.$ Consider the sequence $0\rightarrow L\xrightarrow{f} L\rightarrow T\rightarrow  0$. This gives the following exact sequence
 $$0\rightarrow \Hom_L(\Omega_{L/F},L)\xrightarrow{f}\Hom_L(\Omega_{L/F},L)\rightarrow \Hom_L(\Omega_{L/F},T)$$
 That is, $$0\rightarrow \Der_F(L,L)\xrightarrow{f}\Der_F(L,L)\rightarrow \Der_F(L,T)$$

 By our result in \ref{DerR}, $\Der_F(L)$ is a free $L$-module of $\rank$ equals to $(\dim R+1)$. Therefore, the last exact sequence becomes
 $$0\rightarrow L^{n+1}\xrightarrow{f} L^{n+1}\rightarrow \Der_F(L,T)$$
 Since cokernal of $f$ is $T^{n+1}$, $\rank_T \Der_F(L,T)\geq {n+1}$.
 \end{proof}

     As a corollary of Theorem \ref{BT} we get that $D(R^h)$ is left and right Noetherian ring with global dimension $n$. Moreover, the Bernstein class of $D(R^h)$ is closed under localization at one element of $R^h$.

\section{Proof of theorem \ref{1.5}}
In this section, we present the proof of Theorem \ref{1.5}. Since $R$ is regular of dimension $n$, $E_R(\kappa(P))\cong H^n_P(R)$. We note that $R$ is a $D_F(R)$-module in the Bernstein class, so is $H^n_P(R)$ by \ref{3.12} and this prove first part of \ref{1.5}. We now prove the second part.
\begin{proof}
Let $\kappa(P)=K$, $D=D_F(R)$ and $\widehat D= D_F(\widehat R)$. We note that $E_R(K)\cong H^n_P(R)$ has socle dimension 1. Let the socle be spanned by $t$. Let $\phi$ be the map between $D$ and $E_R(K)$ defined by $\phi(1)=t$. Since $\phi(DP)=0$, so it induces a map $\Bar{\phi}: \frac{D}{DP}\rightarrow E_R(K)$.

Now
\begin{align*}
E_R(K)\otimes_R \widehat R \ & \cong H^n_P(R)\otimes_R \widehat R\\& \cong H^n_{P\widehat R}(\widehat R)\\ &\cong E_{\widehat R}(\frac{\widehat R}{P\widehat R})
\end{align*}

And
\begin{align*}
\frac{D}{DP}\otimes_R \widehat R \ & \cong (D\otimes_R \frac{R}{P})\otimes_R \widehat R\\& \cong (D\otimes_R \widehat R)\otimes_R \frac{R}{P}
\end{align*}
Using \ref{3.5}, $(D\otimes_R \widehat R)\cong \widehat D$. Hence we get the following
\begin{align*}
\frac{D}{DP}\otimes_R \widehat R \ & \cong \widehat D\otimes_R \frac{R}{P}\\ & \cong \widehat D\otimes_{\widehat R} (\widehat R\otimes_R\frac{R}{P})\\ &\cong \widehat D\otimes_{\widehat R} \frac{\widehat R}{P\widehat R}\\ & \cong \frac{\widehat D}{\widehat D\widehat P}
\end{align*}
\ref{3.6} implies $\frac{\widehat D}{\widehat D\widehat P}\cong  E_{\widehat R}(\frac{\widehat R}{P\widehat R})$. Therefore, $E_R(K)\otimes_R \widehat R \cong \frac{D}{DP}\otimes_R \widehat R$. Since $R\rightarrow \widehat R$ is faithfully flat map we get $\frac{D}{DP}\cong E_R(K)$ and this proves our result.
\end{proof}

\section*{Acknowledgement}
The first named author expresses sincere thanks for the invaluable support received from the Government of India through the Prime Minister's Research Fellowship (PMRF ID 1303161), which significantly aided in the completion of this work.
\providecommand{\bysame}{\leavevmode\hbox to3em{\hrulefill}\thinspace}
\providecommand{\MR}{\relax\ifhmode\unskip\space\fi MR }
\providecommand{\MRhref}[2]{
  \href{http://www.ams.org/mathscinet-getitem?mr=#1}{#2}
}

\end{document}